\documentclass[english]{extarticle}
\usepackage[LGR,T1]{fontenc}
\usepackage[latin9]{inputenc}

\makeatletter

\newcommand{\noun}[1]{\textsc{#1}}
\DeclareRobustCommand{\greektext}{%
  \fontencoding{LGR}\selectfont\def\encodingdefault{LGR}}
\DeclareRobustCommand{\textgreek}[1]{\leavevmode{\greektext #1}}
\DeclareFontEncoding{LGR}{}{}
\DeclareTextSymbol{\~}{LGR}{126}

\makeatother

\usepackage{babel}
\begin{document}

\title{A common Misconception about the Categorical Arithmetic}

\author{Giuseppe Raguní}

\date{Universidad Católica de Murcia, Spain - graguni@ucam.edu}
\maketitle
\begin{abstract}
Although the categorical arithmetic is not effectively axiomatizable,
the belief that the incompleteness Theorems can be apply to it is
fairly common. Furthermore, the so-called \emph{essential} (or \emph{inherent})
semantic incompleteness of the second-order Logic that can be deduced
by these same Theorems does not imply the standard semantic incompleteness
that can be derived using the Löwenheim-Skolem or the compactness
Theorem. This state of affairs has its origins in an incorrect and
misinterpreted Gödel's comment at the Königsberg congress of 1930
and has consolidated due to different circumstances. This paper aims
to clear up these questions and proposes an alternative interpretation
for the Gödel's statement.

\end{abstract}

\section{The Categorical Arithmetic}

The categorical arithmetic (\emph{AR}) is a theory where the induction
principle is introduced as a second-order axiom and interpreted according
to the \emph{standard semantics}. This interpretation, briefly called
\emph{full} (since the predicates range over the entire power set
of the universe of discourse), is necessary for the categoricity {[}1{]}. 

\emph{AR} is not a formal theory: it is impossible to dispense with
the meaning of every its formula. This conclusion can be stated as
follow. Since the standard model of \emph{AR} is infinite and unique
under isomorphism, according to the Löwenheim-Skolem Theorem, this
theory cannot be expressed in a semantically complete language\footnote{This also can be concluded by noting that, due to the categoricity,
the compactness Theorem (which applies to every theory with a semantically
complete language), cannot be applicable to \emph{AR}. Here, with
a semantically complete language/axiomatic theory we understand a
system always interpretable if consistent (that is equivalent to affirm
that all the valid formulas are deducible). That is the usual way,
but notice that occasionally others take a quite different meaning
for the same expression (see, \emph{e.g.} {[}2{]}).}. Then, it cannot be formal, because every formal theory has a semantically
complete language {[}3{]}.

\emph{A fortiori}, \emph{AR} is not effectively axiomatizable\footnote{Throughout the paper, we consider as valid the Church-Turing Thesis,
so using always \textquotedbl{}effectively\textquotedbl{} rather than
\textquotedbl{}recursively\textquotedbl{}.}, because every effectively axiomatizable theory is formal.

Since \emph{AR} is finitely axiomatizable, the last conclusion may
surprise. Indeed, the diffuse opinion that every finite set of exhibited
elements always is effectively enumerable, certainly is true when
pure symbols suffice to characterize each one of the elements. But
if we have even one with a meaning, a machine could identify it only
if such meaning is entirely reproducible by mechanical operations
(\emph{i.e.} eliminable) {[}4-5{]}. It is not difficult to realize
that in the case of the \emph{full }second-order induction axiom this
does not occur. Trying to clarify its meaning by a more precise language,
we could resort to that set-theoretical one, representing \emph{AR}
within the Set Theory. Here, introduced a set for the axioms of \emph{AR},
it turns out that the \emph{full} second-order induction principle
is no longer representable by a single axiom but it is equivalent
to an axiomatic scheme capable of generating an uncountable number
of axioms\footnote{Since one axiom for each element of \emph{P(N)} is obtained. Of course,
this is another proof of the non-formality of \emph{AR }{[}6{]}\emph{.}}. 

Therefore, against the widespread view, the incompleteness Theorems
cannot be applied to \emph{AR}. This misconception has dragged on
for too long and has its origins in an incorrect, although absolutely
excusable, Gödel's comment at the Königsberg congress of 1930.

\section{The Gödel's Statement}

According to the editors, the document \emph{{*}1930c} in the third
volume of the \emph{Kurt Gödel\textquoteright s collected works} (1995),
is in all probability the text presented by Gödel at the Königsberg
congress on September 6, 1930 {[}7{]}. In the first part of the document,
Gödel presents his semantic completeness Theorem. After he adds {[}8{]}:
\begin{quotation}
{[}...{]} If the completeness theorem could also be proved for the
higher parts of logic (the extended functional calculus), then it
could be shown in complete generality that syntactical completeness
follows from monomorphicity {[}categoricity{]}; and since we know,
for example, that the Peano axiom system is monomorphic {[}categorical{]},
from that the solvability of every problem of arithmetic and analysis
expressible in Principia mathematica would follow. Such an extension
of the completeness theorem is, however, impossible, as I have recently
proved {[}...{]}. This fact can also be expressed thus: The Peano
axiom system, with the logic of Principia mathematica added as superstructure,
is not syntactically complete.
\end{quotation}
In summary, Gödel affirms that is impossible to generalize the semantic
completeness Theorem to the \textquotedblleft extended functional
calculus\textquotedblright . In fact, in this case also the Peano
axiomatic system, structured with the logic of the \emph{Principia
Mathematica} (\emph{PM}), would be semantically complete. But since
this theory is categorical, it would follow that it is also syntactically
complete. But just this last thing is false, as he \textemdash{} surprise
\textemdash{} announces to have proved.

Now, regardless of what Gödel meant by \textquotedblleft extended
functional calculus\textquotedblright , \emph{this affirmation contains
a} \emph{mistake}. We have in fact two cases:

a) If Gödel understands by \textquotedblleft Peano axiomatic system
structured with the logic of the \emph{PM}\textquotedblright{} an
any type of formal arithmetic theory, the error is precisely to regard
it as categorical.

b) If rather he alludes to the unique categorical arithmetic, that
is \emph{AR}, then Gödel errs applying to it his first incompleteness
Theorem.

Of the two, just the second belief has been consolidating but without
reporting the error. Rather, exalting the merit of having detected
for the first time the semantic incompleteness of the (\emph{full})\emph{
}second-order Logic.

Probably, to forming this opinion has been important the influence
of the following sentence contained in the second edition (1938) of
\emph{Grundzüge theoretischen der Logik} by Hilbert and Ackermann
{[}9{]}:
\begin{quotation}
Let us remark at once that a complete axiom system for the universally
valid formulas of the predicate calculus of second order does not
exist. Rather, as K. Gödel has shown {[}K. Gödel, \emph{Über formal
unentscheidbare Sätze der Principia Mathematica und verwandter systeme},
Mh. Math. Physik Vol. 38 (1931){]}, for any system of primitive formulas
and rules of inference we can find universally valid formulas which
cannot be deduced.
\end{quotation}
The echo of the Gödel\textquoteright s incorrect words at the Congress
pushes the authors (in particular Ackermann, given the age of Hilbert)
to attest that the first incompleteness Theorem concludes directly
the semantic incompleteness of the second-order Logic! False. Furthermore,
as we will try to show, this (true) conclusion really does not appear
to follow by the incompleteness Theorems.

Even in the introductory note of the aforementioned document, Goldfarb
writes {[}10{]}:
\begin{quotation}
Finally, Gödel considers categoricity and syntactic completeness in
the setting of higher-order logics. {[}...{]} Noting then that Peano
Arithmetic is categorical \textemdash{} where by Peano Arithmetic
he means the second-order formulation \textemdash{} Gödel infers that
if higher-order logic is {[}semantically{]} complete, then there will
be a syntactically complete axiom system for Peano Arithmetic. At
this point, he announces his incompleteness theorem: \textquotedblleft The
Peano axiom system, with the logic of Principia mathematica added
as superstructure, is not syntactically complete\textquotedblright .
He uses the result to conclude that there is no (semantically) complete
axiom system for higher-order logic.
\end{quotation}
So interpreting, without the slightest doubt, that Gödel refers to
the second-order categorical arithmetic.

Indeed, today the belief that the incompleteness Theorems can also
apply to \emph{AR} and, above all, that they have as a corollary the
semantic incompleteness of the second-order Logic is widespread. Nevertheless,
it is very rare that someone infers the semantic incompleteness of
the (\emph{full}) second-order Logic in the easy and direct way that
\textemdash{} according to the \emph{b.} interpretation \textemdash{}
Gödel would follow, \emph{i.e.} passing by the (alleged) syntactic
incompleteness of \emph{AR}. Almost all the authors follow the alternative
to prove that the valid formulas of \emph{AR} cannot be effectively
enumerable (see\emph{ e.g.} {[}11{]} and {[}12{]}): by contradiction,
also the true (in the standard model) sentences of the formal (first-order)
arithmetic would be effectively enumerable, against the first incompleteness
Theorem. That is not only more complex but also quite different: the
genuine semantic completeness of a system, simply requires that all
the valid formulas are theorems, not necessarily effectively enumerable
theorems\footnote{This type of semantic incompleteness is called sometimes \emph{essential}
{[}13{]} or \emph{inherent} {[}12{]}, but these adjectives are not
very appropriate because it does not imply the (standard) semantic
incompleteness.}. Actually, the intrinsic non-formality\emph{ }of \emph{AR} entails
that it really makes use of a non-effective deductive method. So,
these proofs really do not conclude the (genuine) semantic incompleteness
of the second-order Logic (as, on the contrary, the use of the Löwenheim-Skolem
Theorem or compactness Theorem can do).

The only explanation of this approach is that the authors are not
sure about the direct applicability of the first incompleteness Theorem
to \emph{AR}. That is not at all surprising in view of the evidences
shown in the previous section; but nothing more is said.

Too respect for the stature of Godel may have affected this state
of affairs, but the main reasons of this misunderstanding are due
probably to ambiguities of the used terminology, both ancient and
modern.

\section{Clearing up the Terms}

The expression \textquotedblleft extended predicate calculus\textquotedblright{}
is for the first time used by Hilbert in the first edition (1928)
of the aforementioned \emph{Grundzüge der theoretischen Logik} where,
with no doubt, indicates the \emph{full} second-order Logic, which
was considered for the first time in the \emph{Principia Mathematica}.
The belief that Gödel, in the aforementioned phrase, refers to the
\emph{AR} theory (explanation \emph{b.}), implies that he, with \textquotedblleft extended
\emph{functional} calculus\textquotedblright , intends the same thing.
But in which work he has shown or at least suggested that the incompleteness
Theorems can apply to the \emph{full} second-order Logic? In none.

In his proof of 1931, Gödel refers to a formal system with a language
that, in addition to the first-order classical logic, allows the use
of non-bound functional variables {[}14{]}. Then he proves that this
is not a real extension of the language, able, in particular, to hinder
the applicability of the semantic completeness Theorem.

In the \emph{1932b} publication, Gödel declares the validity of the
incompleteness Theorems for a formal system (\emph{Z}), based on first-order
logic, with the axioms of Peano and an induction principle defined
by a recursive function. Certainly not a \emph{full} induction. He
adds {[}15{]}:
\begin{quotation}
If we imagine that the system\emph{ Z} is successively enlarged by
the introduction of variables for classes of numbers, classes of classes
of numbers, and so forth, together with the corresponding comprehension
axioms, we obtain a sequence (continuable into the transfinite) of
formal systems that satisfy the assumptions mentioned above {[}...{]}
\end{quotation}
Speaking explicitly of \emph{comprehension axioms}, able to limit
to the countable the number of the sentences, and formal systems.

Finally, in the publication of 1934, which contains the last and definitive
proof of the first incompleteness Theorem, Gödel, having the aim both
to generalize and to simplify the proof, allows the quantification
either on the functional or propositional variables: a declared type
of second-order. However, appropriate comprehension axioms limit again
to infinite countable the number of sentences {[}16{]}. Gödel never
misses an opportunity to point out carefully that always is referring
to a formal system and that the formulas are enumerable {[}17{]}:
\begin{quotation}
Different formal systems are determined according to how many of these
types of variables are used. We shall restrict ourselves to the first
two types; that is, we shall use variables of the three sorts p, q,
r,... {[}propositional variables{]}; x, y, z,... {[}natural numbers
variables{]}; f, g, h, ... {[}functional variables{]}. We assume that
a denumerably infinite number of each are included among the undefined
terms (as may be secured, for example, by the use of letters with
numerical subscripts). {[}...{]} For undefined terms (hence the formulas
and proofs) are countable, and hence a representation of the system
by a system of positive integers can be constructed, as we shall now
do.
\end{quotation}
Therefore, certainly we are not in the \emph{full} second-order. Nevertheless,
in the introduction to the same paper, Kleene, in summarizing the
work of Gödel, does not avoid commenting ambiguously {[}18{]}:
\begin{quotation}
Quantified propositional variables are eliminable in favor of function
quantifiers. Thus the whole system is a form of full second-order
arithmetic (now frequently called the system of \textquotedblleft analysis\textquotedblright ).
\end{quotation}
But he could only mean that the whole system is a \emph{formal} version
(perhaps as large as possible) of the \emph{full} second-order arithmetic.
Maybe is exactly this one the \textquotedblleft extended functional
calculus\textquotedblright{} to which Gödel was referring in the examined
words at the Congress? We will discuss it in the next section. 

Another source of mistake is probably related to use of the term \emph{metamathematics}.
Although Gödel intends it in the modern broad sense that includes
any kind of argument beyond to the coded formal language of Mathematics,
in his theorems always he employs this term limiting it to a formalizable
(though often not yet formalized) use deductive (and, indeed, even
decidable): therefore, only with purpose of brevity.

In the short paper that anticipates his incompleteness Theorems, for
example, Gödel invokes a metamathematics able to decide whether a
formula is an axiom or not {[}19{]}:
\begin{quotation}
{[}...{]} IV. Theorem I {[}first incompleteness Theorem{]} still holds
for all \textgreek{w}-consistent extensions of the system S that are
obtained by the addition of infinitely many axioms, provided the added
class of axioms is decidable, that is, provided for every formula
it is metamathematically decidable whether it is an axiom or not (here
again we suppose that in metamathematics we have at our disposal the
logical devices of PM). Theorems I, III {[}as the IV, but the added
axioms are finite{]}, and IV can be extended also to other formal
systems, for example, to the Zermelo-Fraenkel axiom system of set
theory, provided the systems in question are \textgreek{w}-consistent.
\end{quotation}
But in both the subsequent rigorous proofs, he will \emph{formalize}
this process, which now is called metamathematical, using the recursive
functions, so revealing that, in the words just quoted, he refers
to the usual \textquotedblleft mechanical\textquotedblright{} decidability.
By the same token, even in the theorem that concludes the consistency
of the axiom of choice and of the continuum hypothesis with the other
axioms of the formal Set Theory, he does the same: he uses the metamathematics
only as a simplification, stating explicitly that all \textquotedblleft the
proofs could be formalized\textquotedblright{} and that \textquotedblleft the
general metamathematical considerations could be left out entirely\textquotedblright{}
{[}20{]}.

\section{An Alternative Explanation}

As noted, Gödel has never put in writing that his proofs of incompleteness
may be applied to the uncountable \emph{full} second-order arithmetic
and it looks absolutely not reasonable to believe that he deems it\footnote{I myself have changed my opinion reported in {[}21{]} and {[}22{]}
after a deeper analysis of the \emph{Kurt Gödel\textquoteright s collected
works}.}. In this section, therefore, we will examine the other possibility,
namely the \emph{a.} of the second section. It pretends that Godel
in 1930 believed, mistakenly, categorical a kind of formal arithmetic
and, in consequence of his incompleteness Theorems, semantically incomplete
its language. Is this reasonable (or more reasonable than the previous
case)?

Certainly not for the system considered by Gödel in his first proof
of 1931: in fact, the semantic completeness Theorem applies to it,
as Gödel himself remarks in note n. 55 of the publication {[}23{]}.
Indeed, this is the first time in which the existence of non-standard
models for a formal arithmetic is proved: why Gödel does not report
it? The topic deserves a brief analysis.

More generally than the use of the incompleteness Theorems, the existence
of non-standard models for any formal arithmetic theory can be proved
using the compactness Theorem, the \emph{upward }Löwenheim-Skolem
one or a theorem proved by Skolem in 1933\footnote{ That resolves that the system of the formal (first-order) sentences
that are true in the standard\emph{ }model is non-categorical {[}24{]}.}. The compactness Theorem is due precisely to Gödel (1930) and derives
from his semantic completeness Theorem; but in none of his works Gödel
ever uses it\footnote{This is confirmed by Feferman {[}25{]}}. Moreover,
despite its fundamental importance for the model theory, nobody \textemdash{}
except Maltsev in 1936 and 1941 \textemdash{} uses it before 1945\footnote{This is attested both by Vaught and Fenstad in {[}26{]} and {[}27{]}.}.

Not much more fortunate is the story of the Löwenheim-Skolem Theorem.
The first proof, by Löwenheim (1915), will be simplified by Skolem
in 1920\footnote{The respective works, relative to the first-order classical Logic
(instead of the semantic completeness), are in {[}28{]}.}. In both cases, these theorems are \emph{downward} versions, able
to conclude the non-categoricity of the formal theory of the real
numbers and of the formal Set Theory, but not of the formal Peano
arithmetic. However, Skolem and Von Neumann suspect a much more general
validity of the result {[}29{]}. It seems that also Tarski was interested
to this argument at that time, probably getting the \emph{upward}
version of the Theorem in a seminar of 1928\footnote{This information is due to Maltsev: in {[}30{]} he claims to have
known it by Skolem.}. In any case, the argument continues to have low popularity\footnote{It is significant, for example, that Hilbert did not mention this
theme in the seminary of Hamburg in 1927 {[}31{]}.}, at least until the generalization of Maltsev in 1936 {[}30{]}, which,
including for the first time the \emph{upward} version, will allow
the general conclusion that all the theories equipped with an infinite
model and a semantically complete language are not categorical.

In this context of disinterest for the topic, Gödel not only is no
exception, but his notorious Platonist inclination pushes him to distrust
and/or despise any interpretation that refers to objects foreign to
those that he believes existing independently of the considered theory;
which, in all plausibility, also believes unique. As a matter of fact,
in the introduction of his first paper on the semantic completeness,
he shows to believe categorical even the first-order formal theory
of the real numbers {[}32{]}.

On this basis, one can surmise the following alternative for the option
\emph{a.} When he discovers the non-categoricity of the formal arithmetical
system where his original incompleteness Theorems are applied, Gödel
is not so glad and immediately looks for an extension that, though
formal, is able to ensure the categoricity. Probably he believes to
have identified it in a formal version of the \emph{full} second order
arithmetic: just that one that will be considered in his generalized
proof of the first incompleteness Theorem of 1934 {[}33{]}, where
quantification on the functional and propositional variables are allowed,
while the formality is respected. This hypothesis is consistent with
the fact that Gödel could admit the possibility that this theory uses
a semantically incomplete language, because in both the versions of
his semantic completeness Theorem, he does not allow the use of quantifiers
on functional variables {[}34{]}. Just the planning of this generalization
(literally \emph{extended} to the \emph{functional calculus}) pushes
him, in the meantime, to communicate the result without mentioning
the discovery of the non-standard models. For example, just after
his famous announcement at the Congress of Königsberg on September
6, 1930, Gödel declares {[}35{]}:
\begin{quotation}
(Assuming the consistency of classical mathematics) one can even give
examples of propositions (and in fact of those of the type of Goldbach
or Fermat) that, while contentually true, are unprovable in the formal
system of classical mathematics. Therefore, if one adjoins the negation
of such a proposition to the axioms of classical mathematics, one
obtains a consistent system in which a contentually false proposition
is provable.
\end{quotation}
Indubitably here, admitted the soundness, he has in mind a formal
theory that, being syntactically incomplete and categorical, has a
semantically incomplete language. 

The first surprise comes to him with the Skolem\textquoteright s proof
of 1933: neither the system of the formal (first-order) sentences
that are true in the standard\emph{ }model is categorical {[}24{]}.
A first disturbing clue that the non-categoricity covers all the formal
(with at least an infinite model) systems, regardless of the syntactic
completeness or incompleteness. Gödel, in reviewing the Skolem\textquoteright s
paper, laconically observes \textemdash{} finally! \textemdash{} that
a consequence of this result, that is the non-categoricity of the
formal Peano arithmetic, was already derivable from his incompleteness
Theorems {[}36{]}. Later, in any work (not only in the cited generalization
of 1934), he always will ignore the issue of the categoricity, nor
ever will return to state that by his incompleteness Theorems can
be derived the semantic incompleteness of some language or theory.

Ultimately, the Henkin\textquoteright s Theorem of 1950 {[}3{]} will
prove that in every formal system (and so, anywhere the incompleteness
Theorems could be applied) there is semantic completeness of the language
and therefore, if at least an infinite model exists, there cannot
be categoricity.

\section{Conclusions}

Since the categorical arithmetic is not effectively axiomatizable
and any type of formal arithmetic is not categorical, the text of
the Gödel\textquoteright s communication at the conference in Königsberg
on September 6, 1930 (never published by him) contains a mistake.
In the common understanding this error is not reported and thus it
is wrongly believed that: a) the incompleteness Theorems also can
be applied to the categorical arithmetic; b) the semantic incompleteness
of the second-order Logic is a consequence of the incompleteness Theorems. 

The previous interpretation is untenable, nor supported by the Gödel's
publications. As a matter of fact, the semantic incompleteness of
the second-order Logic is due to the fact that this language allows
to get the categoricity of theories (not only \emph{AR}) equipped
with at least an infinite model. And this is independent of the syntactic
incompleteness of the formal arithmetic. By the incompleteness Theorems
it is only possible to derive the so-called \emph{essential} (or \emph{inherent})
semantic incompleteness of the second-order Logic, which however does
not imply the standard semantic incompleteness.

As an alternative interpretation of the manuscript in question, it
is very plausible that Gödel was referring to a formal arithmetic
(later specified in his proof of 1934) in which the quantification
on the functional and propositional variables is allowed. If so, in
1930 he believed that this theory was categorical and, as a consequence
of its syntactic incompleteness, equipped with a semantically incomplete
language. This explanation is consistent with the fact that both the
version of his semantic completeness Theorem cannot be applied to
this system, due to the said quantification.

We wish to emphasize that this alternative in no way shades the luster
of Gödel, because it makes no sense to pretend that he, in 1930, could
know that every formal system, equipped with at least one infinite
model, is not categorical. Conversely, it absolves him from a blunder
and also explains why, becoming more and more evident, as time passes,
the difficulty for the condition of categoricity, he never will repeat
alike affirmations. On the other hand, Gödel never corrected the phrase
presumably because he was not worrying about rectifying an unpublished
text.

Finally, about the possible syntactic (and, by the categoricity, also
semantic) completeness of the categorical arithmetic, we just observe
that it would not be incompatible with the fact that the \emph{language}
of this theory is semantically incomplete. In fact, although an axiomatic
system that uses a semantically complete language always is semantically
complete, the reverse is not always true {[}37{]}.

\end{document}